\documentclass[draft]{article}
\usepackage[leqno,intlimits]{amsmath}
\usepackage{amssymb}
\usepackage{amsthm}
\usepackage{hyperref}

\newtheorem{tm}{Theorem}[section]
\newtheorem{lm}[tm]{Lemma}
\newtheorem{pr}[tm]{Proposition}
\newtheorem{cor}[tm]{Corollary}

\newtheorem{Def}[tm]{Definition}

\numberwithin{equation}{section}

\newcommand\kon{r} 
\newcommand\mom{m} 
\newcommand*{\un}[1]{\underline{#1}}
\newcommand*{\Zb}{\mathbb Z}
\newcommand*{\Rb}{\mathbb R}
\newcommand*{\om}{\omega}
\newcommand*{\ze}{\zeta}
\newcommand*{\de}{\delta}
\newcommand*{\la}{\lambda}
\newcommand*{\ba}{\begin{aligned}}
\newcommand*{\ea}{\end{aligned}}
\newcommand*{\be}{\begin{equation}}
\newcommand*{\ee}{\end{equation}}
\newcommand*{\e}[1]{\text{\rm e}^{#1}}
\newcommand*{\vr}{\varrho}
\newcommand*{\Ev}{{\bf E}}
\newcommand*{\Pv}{{\bf P}}
\newcommand*{\Vv}{{\text{\bf Var}}}
\newcommand*{\Oc}{\mathcal O}
\newcommand*{\lc}{\lceil}
\newcommand*{\rc}{\rceil}
\newcommand*{\lf}{\lfloor}
\newcommand*{\rf}{\rfloor}
\newcommand*{\di}{\,\text{\rm d}}
\newcommand*{\al}{\alpha}
\newcommand*{\ga}{\gamma}
\newcommand*{\wt}{\widetilde}
\newcommand*{\hop}{\bigskip\noindent}
\newcommand*{\hip}{\smallskip\noindent}

\begin{document}

\title{Order of current variance  and diffusivity  in the rate one totally asymmetric zero range process}
\author{
M\'arton Bal\'azs\thanks{MTA-BME Stochastics Research Group},
J\'ulia Komj\'athy \thanks{Budapest University of Technology and
Economics \newline Bal\'azs was partially supported by the Hungarian
Scientific Research Fund (OTKA) grants K60708, TS49835, and F67729,
the Bolyai Scholarship of the Hungarian Academy of Sciences.}}
\maketitle
\begin{abstract}
We prove that the variance of the current across a characteristic is
of order $t^{2/3}$ in a stationary constant rate totally asymmetric
zero range process, and that the diffusivity has order $t^{1/3}$.
This is a step towards proving universality of this scaling behavior
in the class of one-dimensional interacting systems with one conserved 
quantity and concave hydrodynamic flux.
The proof proceeds via couplings to show the corresponding moment
bounds for a second class particle. We build on the methods
developed in \cite{se2/3} for simple exclusion. However, some
modifications were needed to handle the larger state space. Our
results translate into $t^{2/3}$-order of variance of the tagged
particle on the characteristics of totally asymmetric simple
exclusion.
\end{abstract}

\noindent {\bf Keywords:} Constant rate totally asymmetric zero
range process, diffusivity, current fluctuations, second class
particle

\hop
{\bf 2000 Mathematics Subject Classification:} 60K35, 82C22

\section{Introduction}

The constant rate totally asymmetric zero range process (TAZRP) is a
Mar\-kov process that describes the motion of particles in the one
dimensional integer lattice $\Zb$. Namely, if any particles are
present at a site, then one of them jumps one site to the right with
rate $1$ independently of particles at other sites. Later we will
give a construction of TAZRP in terms of Poisson clocks. This
process is among the interacting particle systems introduced in
\cite{spi}. Our paper examines the net particle current seen by an
observer moving at the characteristic speed of the process. The
characteristic speed is the speed at which perturbations travel in
the system and can be determined e.g.\ via the hydrodynamic limit.
The process is assumed to be in one of its extremal stationary
distributions, i.e.\  with independent geometrically distributed
occupation variables at each site with parameter $\frac{1}{1+\vr}$,
implying density $\vr>0$ of particles. We prove that the net
particle current across the characteristic, that counts the number
of particles that pass the observer from left to right minus the
number that pass from right to left during time interval $(0,t]$,
has variance of order $t^{2/3}$.

The seminal papers of Baik, Deift and Johansson \cite{bdj}  and
Johansson \cite{1/3}  gave the first rigorous proofs of such
fluctuations in the last-passage version of Hammersley's process and
in the totally asymmetric simple exclusion process (TASEP) started
from jam-type initial configuration. The correct order was verified
to be $t^{1/3}$, and the limiting fluctuations were found to obey
Tracy-Widom distributions from
 random matrix theory. Later a last-passage representation was also found for a stationary
TASEP \cite{spohn}, and then the Tracy-Widom limit  proved for the
current across the characteristic in that setting \cite{ferspohn}.

 Cator and Groeneboom  \cite{cuberoot} used a probabilistic
approach to prove the correct order of these fluctuations for
Hammersley's process, Bal\'azs, Cator and Sepp\"al\"ainen
\cite{third} obtained the same result for the last passage
representation of the stationary TASEP. Then Bal\'azs and
Sepp\"al\"ainen in \cite{se2/3} proved that the current variance is
of order $t^{2/3}$ in the (non-totally) asymmetric simple exclusion
process (ASEP).

There are generalizations of the result for the ASEP for those
processes in which particles can jump more than one site, most
recently by Quastel and Valk\'o \cite{quava, quava2}.

Physical reasoning, e.g.\ \cite{vBKS85}, suggested universality in the sense that similar scaling of
the current across the characteristic
 should occur in other systems as well with one conserved quantity.

  To our knowledge, the present
article is the first one that generalizes the result of \cite{se2/3}
towards another system in which more than one particle per site is
allowed. Due to the connection between the TASEP tagged particle and
the TAZRP particle current, our results also apply to tagged
particle fluctuations of the TASEP. As in \cite{se2/3}, our
arguments are entirely probabilistic and utilize couplings of
several processes and bounds on second class particles. Informally
speaking, {\sl second class particles} are perturbations in the
system that do not disturb the motion of the regular particles but
are influenced by the ambient system. Precise definition will be
given after the construction of the coupled processes in Section
\ref{sc:prore}. As the state space is larger than in the ASEP, some
delicate modifications are required when we couple several processes
together.

Fluctuation results for asymmetric exclusion processes have also
been stated in terms of  a quantity called the {\sl diffusivity}
$D(t)$. The link between current variance and diffusivity can be
described also using the concept of the
 second class particle:  the variance of
the current is the expected absolute deviation of the second class
particle, while $tD(t)$ is the variance of the second class
particle.  For the TAZRP we also obtain the correct order $t^{1/3}$
for the diffusivity.

When the observer's speed $V$ is different from the characteristic
speed, the fluctuations are normal, i.e.\ the current variance is of
order $t$. This was proved by Ferrari and Fontes in \cite{se} for
the ASEP and by Bal\'azs for some other processes in \cite{fluct}.
As a consequence we also reproduce these results concerning the
normal fluctuations for TAZRP.

This paper is the continuation of \cite{se2/3} and might give some
ideas that could lead to the treatment of some other models. The
method seems to be robust enough to be able to generalize it to
other, even more complicated processes. Proceeding in this direction
towards universality is subject to future work.

\section{The constant-rate zero range process and the results}\label{sc:prore}


\noindent{\bf Construction of the process and the second class
particles}

\hip The constant-rate totally asymmetric zero range process (TAZRP)
is a Markov process on the state space
$\Omega=\{0,1,2,\dots\}^\Zb={\Zb^+}^\Zb$. Given a state
$\un\om=\{\om_i\}_{i\in\Zb}\in\Omega$, the following jumps can
happen independently at different sites:
\be
(\om_i,\,\om_{i+1})\longrightarrow(\om_i-1,\,\om_{i+1}+1)\text{ with rate }
\mathbb I_{(\om_i>0)}.\label{eq:up}
\ee

We interpret the process as representing unlabeled particles that
execute nearest-neighbor random walks on $\Zb$, in such a way that
at a site where at least $1$ particle is present only the topmost
one is allowed to jump to the right with rate $1$. The value
$\om_i(t)=0$ means that site $i$ is vacant at time $t$. The state of
the entire process at time $t$ is then $\un\om(t)$.

A rigorous construction of this process is done by giving each site
$i$ a rate $1$ Poisson process $N_{i\to i+1}$ on the time line
$[0,\infty)$. The processes $\{ N_{i\to i+1}, i\in\Zb\}$ are
mutually independent, and also independent of the initial
configuration $\un\om(0)$.  The rule of evolution is that when
$N_{i\to i+1}$ jumps, the topmost particle is moved from $i$ to
$i+1$ if $i$ is occupied by at least one particle.  Thus the rates
\eqref{eq:up} are realized.

Let $\mu_\vr$ denote the geometric measure with mean $\vr$, i.e.\
 $\mu_\vr\{k\}=\Big(\frac{\vr}{\vr+1}\Big)^k
\frac{1}{\vr+1}$ on the set $\Zb^+$, and let
$\un\mu_\vr=\mu_\vr^{\otimes\Zb}$ be the i.i.d.~geometric product
measure with marginals $\mu_\vr$ on $\Omega$. It is known that the
measures $\{\un\mu_\vr\,:\,0\leq\vr\}$ are the extreme points of the
convex set of invariant distributions for the process that are also
invariant under spatial translations.

It is convenient to embed the zero range process in a height process
that represents a wall of adjacent columns of bricks. On top of each
interval $[i,i+1]$ sits a  column of bricks with
 height $h_i\in\Zb$.  The entire height configuration
is $\un h=\{h_i\}_{i\in\Zb}$,  restricted to satisfy \be 0\leq
h_{i-1}-h_i \quad\text{ for each $i$} \label{eq:h} \ee so that the
wall slopes downward to the right. Let the Poisson processes govern
the evolution of the heights:
 when $N_{i\to i+1}$ jumps add a brick
on top of the column on $[i,i+1]$. But suppress every step that
leads to a violation of \eqref{eq:h}.

Given an initial particle configuration $\un\om(0)$,
define an initial height configuration by
\[
h_i(0)=\begin{cases}
\sum_{j=i+1}^0\om_i(0)&\text{for }i<0,\\
0&\text{for }i=0,\\
-\sum_{j=1}^i\om_i(0)&\text{for }i>0.
\end{cases}\]
Let the heights evolve, and define
\[
\om_i(t)=h_{i-1}(t)-h_{i}(t).
\]
Then this process $\un\om(t)$ is exactly the TAZRP constructed
earlier, and the height increment $h_i(t)-h_i(0)$ is the net
particle current across the bond $(i,i+1)$.

Below we run several processes started from different initial
configurations together governed by the
{\sl basic coupling} in which we use common Poisson clocks for them.
The first observation is that this coupling preserves monotonicity
among both particle and height configurations. Ordering is defined
sitewise: for particle configurations
 $\un\eta\leq\un\om$ means that  $\eta_i\leq\om_i$  for each $i\in\Zb$,
and similarly for height configurations
 $\un g\leq\un h$ if  $g_i\leq h_i$    for each $i\in\Zb$.
The basic coupling has the following property, called attractivity:
\[
\un\eta(0)\leq\un\om(0)\Longrightarrow\un\eta(t)\leq\un\om(t)\quad\text{and}
\quad\un g(0)\leq\un h(0)\Longrightarrow\un g(t)\leq\un h(t)
\]
for all $t>0$.

We use the following terminology:  if we have two coupled zero range
processes  $\un\eta(t)\leq \un\om(t)$, then there are $\om_i-\eta_i$
pieces of {\sl $\om-\eta$ second class particles} at each site $i$.
The joint process $(\un\eta(\cdot),\un\om(\cdot))$ can be
constructed from a two-class process: (i) The first class particles
 $\un\eta$ obey the TAZRP dynamics as described earlier.
(ii) The second class particles also obey the Poisson clocks when
they can, but they are only allowed to jump if there are no more
first class particles at site $i$, i.e.\ their jump rate at site $i$
is $\mathbb I_{(\eta_i=0)}\cdot \mathbb I_{(\om_i\ge 1)}$ and only
the topmost second class particle is allowed to jump. We emphasize
here a different behavior of the zero range second class particles
and the ones in the simple exclusion process, i.e.\  a jump of a
first class particle to site $i$ in TAZRP does not affect the number
of second class particles at site $i$. In the ASEP case, where at
most one particle is allowed at one site, they interchange sites,
implying a left jump of the second class particle even in the
totally asymmetric case. Hence the TAZRP model second class
particles are allowed to jump only to the right.

Let $\un\de_i\in\Omega$ denote a configuration that has only a
single particle  at site $i$. For $\un\eta\in\Omega$ we can
legitimately define
 $\un\eta^+=\un\eta+\un\de_0$. In this situation
we say that there is  a single second class particle between
$\un\eta^+$ and $\un\eta$ at site 0. Since the basic coupling
conserves it, there is always a site $Q(t)$ such that
\[
\un\eta^+(t)=\un\eta(t)+\un\de_{Q(t)}.
\]
 $Q(t)$ is the position of the second class particle at time $t$,
which performs a nearest neighbor walk, influenced by the ambient process
$\un\eta(\cdot)$.

It is convenient to also
have the notion of a second class \emph{anti}particle at position
$Q_a(t)$ in a process $\un\om(t)$.  This means  that
$Q_a(t)$ is the location of the single discrepancy
between two processes $\un\om(t)$ and  $\un\om^-(t)$ that are
 started so that
$\un\om^-(0)=\un\om(0)-\un\de_i\ge 0$ where $i=Q_a(0)$. In the
proofs we will couple more than two processes and this flexibility
will be convenient.


\hop
{\bf Current fluctuations and diffusivity}

\hip
Let $[x]$ denote the first integer from $x$ towards the origin,
in other words  $[x]=\lf x\rf$ (floor) when $x\geq0$
 and $[x]=\lc x\rc$ (ceiling)
when $x<0$.
For a speed value $V\in\Rb$  define \be
J^{(V)}(t)=h_{[Vt]}(t),\label{eq:jvdef} \ee the height of the column
over interval $[\,[Vt],[Vt]+1]$ at time $t$. The normalization
$h_0(0)=0$ implies that
  $J^{(V)}(t)$ is the total net particle current seen by an observer moving
at speed $V$ during time interval $[0,t]$.
One can compute the density $\vr$ stationary expectation
\be\label{eq:jvt} \Ev^\vr J^{(V)}(t)= \frac{\vr}{\vr+1} t -\vr [Vt] \ee
by writing  a martingale for
 $h_0(t)$ and then adding in
$h_{[Vt]}(t)-h_0(t)$ which counts particles between sites $0$ and
$[Vt]$.

Our results are based on an interplay between currents and second
class particles.  One key fact  is the coming connection. To
simplify notation we introduce a new measure
\[\hat\mu_\vr(k)=\frac{1}{\Vv^\vr(\om_0)} \sum_{y=k+1}^{\infty}
(y-\vr) \mu_\vr(y).\] It can be checked easily that $\hat \mu$ is
again a probability measure for any probability measure $\mu$. It
arises from \cite[Theorem 2.2]{varj2nd}. An easy computation gives
that in the case of TAZRP
 $\hat
\mu_\vr(k)=(k+1)\Big(\frac{\vr}{1+\vr}\Big)^k
\Big(\frac{1}{1+\vr}\Big)^2$, which is the distribution of the sum
of two independent geometrically distributed random variables with
mean $\vr$.  We will use the configuration of a coupled pair of
processes started  in a product measure where all sites have
independent initial marginals $\mu_\vr$ except the origin where the
initial distributions are $\hat \mu_\vr$ and $\hat \mu_\vr+1$
respectively, independently of other sites. Denote this coupling
measure of the two states $\un \om^-(0)$ and $\un \om(0)$ by $\un
{\hat \mu}_\vr$.

From now on $\Ev^\vr$ and $\Vv^\vr$ denotes expectation and variance
in the stationary process with density $\vr$ (process started from
$\mu_\vr^{\otimes \Zb}$) while $ \Ev$ and $ \Vv $ denotes
expectation and variance in a coupled pair of processes started from
$\un{\hat \mu}_\vr$.
\begin{pr}\label{pr:se}
Let $\un\om(\cdot)$ and $\un \om^-(\cdot)$ be a coupled pair of
TAZRP started from the distribution $\un{\hat \mu}_\vr$. Let
$Q_a(\cdot)$ be a second class antiparticle between the processes
starting at $Q_a(0)=0$. Let the model evolve according to the basic
coupling. Then the variance of the current for any $V\in \Rb$ in the
stationary process can be rewritten as:
 \be
\begin{split}
\Vv^\vr(J^{(V)}(t))&=\Vv^\vr(\om_0) \Ev\bigl(\,|\,[Vt]-Q_a(t)|\,\bigr)\\
&=\vr (1+\vr)\Ev\bigl(\,|\,[Vt]-Q_a(t)|\,\bigr).
  \end{split}
\label{eq:varjse} \ee Also,
 \be \Ev\bigl(Q_a(t)\,)=V^\vr t =\frac{1}{(1+\vr)^2}t
\label{eq:Qmean}\ee holds.
\end{pr}
 In
article  \cite{varj2nd} identities \eqref{eq:varjse} and
\eqref{eq:Qmean} are proved for a class of attractive models that include the asymmetric simple
exclusion,
 zero-range and bricklayer processes.

The interesting current fluctuations occur at the
\emph{characteristic speed} $V^\vr=\frac{1}{(1+\vr)^2}$. From
\eqref{eq:Qmean} we see that this is the average speed of the second
class particle. The characteristic speed also appears in the PDE
theory of the conservation law
 \be\vr_t+f(\vr)_x=0 \label{eq:conslaw}\ee that is obtained via the
Eulerian hydrodynamic limit of TAZRP. Here the hydrodynamic flux
takes the form
 \be \label{eq:flux}f(\vr)=\Ev^{\vr}\big(\mathbb I_{(\om>0)}\big)=\frac{\vr}{1+\vr}.
 \ee
 At constant density $\vr$ the
characteristic speed is $f'(\vr)=V^\vr$ at which perturbations of
entropy solutions of \eqref{eq:conslaw} travel.
For more information we refer  to \cite{cl} for
 hydrodynamic limits
of interacting particle systems.

Let us return to the stationary TAZRP with
Geometric$\big(\frac{1}{1+\vr}\big)$ occupation variables at each
fixed time. A basic object for understanding space-time correlations
 is the  {\sl two point function}
$S(i,t)=\Ev^\vr\big(\om_i(t)\om_0(0)\big)-\vr^2$. Due to
\cite{varj2nd} it can be written as the transition probability of
the second class particle: \be\label{eq:twopoint} S(i,t) =
\Ev^\vr\Big(\mathbb I\{Q(t)=n\} \cdot \sum_{z=\om_0+1}^{\infty}
(z-\vr) \frac{\mu_\vr(z)}{\mu_\vr(\om_0)} \Big) =
\Vv^\vr(\om_0)\cdot \Pv\big\{Q(t)=n \big\}. \ee
The sum under the expectation is regarded as a Radon-Nikodym
derivative that takes our second class particle in the initial
setting $\un{\hat\mu}_\vr$, where the second class particle resides
initially at the origin. Note that use of the term ``transition
probability'' is not meant to suggest that $Q(t)$ is a Markov
process.

From \eqref{eq:twopoint} and \eqref{eq:Qmean} follows
\[
 \sum_{i\in\Zb}i S(i,t) = \Vv^\vr(\om_0) V^\vr t=\vr(1+\vr) V^\vr t.
\]
 The  {\sl diffusivity}
is by definition a normalized second moment of the two-point
function:
\[
D(t)=\frac1{t\vr(1+\vr)} \sum_{i\in\Zb} (i-V^\vr t)^2
S(i,t).
\]
Consequently  the  diffusivity can also
be expressed in terms of the variance of the second class particle:
\[
D(t)=t^{-1}\Vv(Q(t)\,).
\]
Below we show that, similarly to the ASEP, the TAZRP second class
particle is also {\sl superdiffusive} with variance of order
$t^{4/3}$.

We can now state the main theorem,
 a moment bound for the second class
particle in a pair started from $\un{\hat\mu}_\vr$.

\begin{tm}\label{tm:main}
For the constant rate TAZRP $\un\om(\cdot)$  started in $\un{\hat
\mu}_\vr$ distribution for any $\vr\ge0$ there exist constants
$0<t_0, C<\infty$ such that for real $1\le \mom<3$ and $t\ge t_0$
\be C^{-1}\;\le\; \Ev\biggl\{\, \Bigl\lvert\,\frac{Q_a(t)-V^\vr
t}{t^{2/3}}\Bigl\lvert^m\,\biggr\} \;\le\; C. \label{eq:mainbd}\ee
\end{tm}

 Let us
 simplify notation to $J^\vr(t)=J^{(V^\vr)}(t)$ for the
current across the characteristic.
Via \eqref{eq:varjse},
taking $\mom=1$ gives the variance of this current.

\begin{cor}
 With assumptions as in Theorem \ref{tm:main}, for $t\ge t_0$
 the variance
of the current across the characteristic satisfies
\[
C^{-1}t^{2/3} \le {\Vv^\vr(J^\vr(t))} \le C t^{2/3}.
\]
\end{cor}

Taking $\mom=2$
 identifies the order of the diffusivity.

\begin{cor} With assumptions as in Theorem \ref{tm:main},
for $t\ge t_0$
\[
C^{-1}t^{1/3} \le D(t) \le C t^{1/3}.
\]
\end{cor}

The upper bound of \eqref{eq:mainbd} is not valid
for all $t>0$. For small $t$ the variance of $Q(t)$
is of order $t$ because the second class particle
 is likely to have experienced at most
one  jump.

The above results are completely analogous to those in \cite{se2/3}.

Next, another important consequence of the theorem is the Weak Law
of Large Numbers for $Q(t)$. Markov's inequality implies
\[
\Pv\Big(\Big|\frac{Q(t)}{t} -V^\vr \Big|>\varepsilon\Big)\le
\frac{\Ev|Q(t)-[V^\vr t]|+1}{\varepsilon t}\le \frac{C
t^{2/3}}{\varepsilon t}\to 0.
\]
We mention here that the moment
bounds of \eqref{eq:mainbd} together with the fact that $Q(t)$ is
dominated by a Poisson process are strong enough to derive the
Strong Law of Large Numbers as well.

Proposition \ref{pr:se} and Theorem \ref{tm:main} imply that the
current variance is of order $t$ for all $V\neq V^\vr$:
\[
\ba
\lim_{t\to \infty} \frac{\Vv^\vr J^{(V)}(t)}{t}&= \lim_{t\to\infty}
\frac{1}{t} \Vv^\vr(\om_0) \Ev(|Q(t)- [Vt]|) \\&= \lim_{t\to \infty}
\frac{1}{t}\Vv^\vr(\om_0) \Big(\Ev|Q(t)-[V^\vr t] + [V^\vr
t]-[Vt]|\Big)\\&= \Vv^\vr(\om_0)|V^\vr-V|.
\ea
\]

The Central Limit Theorem for $J^{(V)}(t)$ also holds in the
following form:
\[
\Pv\Big(\frac{J^{(V)}(t) - \Ev^\vr
J^{(V)}(t)}{\sqrt{t \Vv^\vr(\om_0)|V^\vr-V|}}<x\Big)\rightarrow
\Phi(x)\mbox{ for all } x \in \mathbb R, \quad V\neq V^\vr,
\]
where $\Phi(x)$ denotes the standard normal distribution.

 A more important consequence of Theorem
\ref{tm:main} concerns the fluctuations of the tagged particle  of
the TASEP near the characteristics. To be able to explain the result
we introduce some new notation: first, we consider a stationary
TASEP with independent Bernoulli($\alpha$) occupation variables at
each site which we condition on the event that the origin is
occupied. Next we label the particles in an increasing order such a
way that the one at the origin becomes label $0$. We denote the site
of the $i$-th particle at time $t$ by $R^{E}_i(t)$. Then we let the
system evolve according to the TASEP evolution with particles
jumping to the left and not to the right. Rates of these jumps are
$\mathbb I_{\{\om^E_i=1\}} \cdot \mathbb I_{\{\om^E_{i-1}=0\}}$,
i.e.\  the particles are allowed to jump one to the left if the
neighboring site on the left is empty. The particles in the TASEP
cannot pass each other, hence the ordering is preserved over time.
As we conditioned on the origin being occupied, $R^{E}_0(0)=0$ and
we choose the usual height representation of the process with
$h^{E}_0(0)=0$ and the height sloping to the left. Since the
particles jump to the left, this height function is decreasing in
time.

Second, we look at the process of the differences, namely the empty
sites between two neighboring particles in the exclusion process.
Define $\om_i(t):=R^E_i(t)-R^E_{i-1}(t)-1$. Then the process
$\un\om(t)$ is a stationary TAZRP with parameter $\vr=1/\alpha-1$,
where particles jump to the right. The height representation $\un
h(t)$ of the zero range process can be translated to the TASEP,
namely \be\label{contact} R^{E}_k(t)=-h_k(t)+k. \ee For $k=[V^\vr
t]$, we get $R^{E}_{[V^\vr t]}(t) = -J^{\vr}(t)+[V^\vr t]$. Theorem
\ref{tm:main} leads to \[ C^{-1}\le\liminf_{t\to \infty}
\frac{\Vv\big(R^{E}_{[V^\vr t]}(t)\big)}{t^{2/3}} \le \limsup_{t\to
\infty} \frac{\Vv\big(R^{E}_{[V^\vr t]}(t)\big)}{t^{2/3}} \le C \]

To give an interpretation of this result, we recall that the
characteristic speed in the TAZRP equals $V^\vr=
\frac{1}{(1+\vr)^2}= \alpha^2,$ so equation \eqref{contact} at site
$[V^\vr t]$ turns into $R_{[\alpha^2 t]}(t)=- J^{\vr}(t) + [\alpha^2
t].$ Since the mean distance between particles in the TASEP is
$1/\al$, the $[\alpha^2 t]^{\text{th}}$ particle started at time $0$
somewhere near $\alpha t$. Burke's Theorem implies that each
particle marginally moves to the left following a Poisson process
with parameter $1-\al$, hence this particle should be at time $t$
somewhere near $(2\al-1) t$, which is exactly the characteristic
position of the exclusion process. So, to put the result
$\Vv(R^{E}_{[\alpha^2 t]}) \sim t^{2/3}$ in words we can say that
the position of the tagged particle that is expected to be on the
characteristics has fluctuation of order $t^{1/3}$.

\hip
 The rest of the paper is devoted to the proof of Theorem
\ref{tm:main}, first the upper bound and then the lower bound.

\section{Upper bound}
\label{sc:ub} Proof of the upper bound utilizes couplings of several
processes. Before deriving the upper bound we introduce some
preliminaries on the couplings.


\hop
{\bf Couplings}

\hip
We start by describing an initial distribution of three coupled
processes and with labels attached to second class particles between
two of the three processes.
\begin{Def} For two probability measures $\mu_1\le \mu_2$ denotes stochastic
domination. \end{Def} \noindent Equivalent to stochastic domination
are:
\begin{itemize}
\item For the cumulative distribution functions $F_1(x) \ge F_2(x)$
holds for all $x\in \mathbb R$.
\item  There exists a coupling measure
$\mu(x,y)$ with marginals $\mu_1$ and $\mu_2$
such that $\mu\big\{(x,y)\,:\, x \le y\big\}=1$.
\end{itemize}

 Fix densities $0<\la<\vr$. In the next sections we want to
 couple models with different densities and distribution so here we prove stochastic domination of the different measures.
 \begin{lm} The measures $\mu_\vr$ and $\hat \mu_\vr$ are monotone in $\vr$, i.e.
 \begin{enumerate}
 \item  $\mu_\la \le \mu_\vr$ for all $\la \le\vr$,
  \item $\hat\mu_\la \le \hat \mu_\vr$ for all $\la \le\vr$
\item Furthermore, $\hat \mu_\vr \ge \mu_\vr$ for all $\vr>0$.
\end{enumerate}
 \end{lm}
 \begin{proof}
 The cumulative distribution function of the geometric distribution with mean $\vr$ is
 \[F_\vr(z)=\sum_{k=0}^{z}\mu_\vr(k)=\sum_{k=0}^{z}\frac{\vr^k}{(1+\vr)^{k+1}}=
 1-\Big(\frac{\vr}{1+\vr}\Big)^{z+1}.\]
Differentiation with respect to $\vr$ gives that the function is
monotone decreasing in $\vr$, meaning that $\la\le \vr$ implies
$F_\vr(z)\le F_\la(z)$. Next, $\hat \mu_\vr$ is the distribution of
the sum of two independent, geometrically distributed random
variables, implying $\widehat F_\vr(z)\le \widehat F_\la(z)$. The
third statement comes from the fact that a random variable of
distribution $\hat \mu_\vr$ is the sum of two independent random
variables of distribution $\mu_\vr$.
\end{proof}
\begin{Def}\label{def:measures}
\hip
\begin{itemize}
\item $\nu$ is the coupling measure with marginals $\mu_\la$ and
$\mu_\vr$ with $\nu\{(\eta, \om): \eta\le \om\}=1$. Let
$\un\nu:=\nu^{\otimes \Zb}$ be the product measure with marginals
$\nu$ for each site.
\item $\hat \nu$ is the coupling measure with marginals $\hat
\mu_\la$ and $\hat \mu_\vr$ with $\hat \nu\{(\eta,\om): \eta\le
\om\}=1$.
\end{itemize}
\end{Def}
We now describe the initial distribution. Let $(\un\eta, \un\om^-)$
be distributed in the product measure of marginals $\nu$ for each
site $i\neq 0$ and $\hat \nu$ at the origin. We also couple a third
process denoted by $\un\om$ to these two processes with
$\om(0)=\om^-(0)+\delta_0$ i.e. we add an extra second class
particle to the process $\un\om^-$. Saying differently, we have a
second class antiparticle on $\un\om$. Let $\bar {\un\nu}$ denote
the resulting joint distribution of $(\un\eta, \un\om^-, \un\om)$.
Then $\bar{\un\nu}$-a.s.\ we have the following initial conditions
for the processes $\un\eta$, $\un\om^-$ and $\un \om$ :
\begin{itemize}
\item $\un\eta\leq\un\om^-\leq \un \om$,
\item There is a second class particle between $\un\om$ and $\un\om^-$ at
the origin at time $0$, implying that there is at least one second
class particle between $\un\om$ and $\un \eta$ at the origin.
\item  There are infinitely many $\om-\eta$ second class particles on both sides of the origin.
\end{itemize}

 Label the  $\om-\eta$ second class particles with integers in
an increasing fashion from left to right, giving label 0 to the
topmost second class particle initially at the origin. $X_k(0)$
denotes the position of second class particle with label $k$, so
that \be\dotsm \le X_{-1}(0)\le X_0(0)=0<X_{1}(0)\le\dotsm
\label{eq:Xlabels}\ee

Let these configurations evolve from the initial distribution $\un
{\bar\nu}$
 with common Poisson clocks.
The rule for the labeling is that whenever  a second class particle
jumps from site $i$, it is always the highest labeled one at $i$ and
it arrives at the lowest place at site $i+1$. This way the ordering
\eqref{eq:Xlabels} is kept for all time, and the basic coupling
preserves the ordering $\un\eta(t)\leq\un\om^-(t) \leq\un\om(t)$
a.s.\ for all times $t$. The topmost $\om-\eta$ second class
particle that starts at the origin is of special importance to us,
so we set $X(t)=X_0(t)$. The effect of the coupling is that the
$\om-\om^-$ particle (denote its position from now on by $Q_a(t)$)
cannot pass $X(t)$. Precisely saying, the ordering $Q_a(t) \le X(t)$
is preserved over time. The argument for this fact is that initially
$Q_a(0) \le X(0)$ holds. Since both particles are only allowed to
jump to the right $Q_a$ can pass $X$ if and only if they are at the
same site. In this case $X(t)$ jumps if $\eta_{X(t)}(t)=0$ but
$Q_a(t)$ jumps only if $\om^-_{Q_a(t)}=0$, and the coupling
$\un\eta\le \un\om^-$ guarantees that the second event is part of
the first one, i.e. the jump of $Q_a(t)$ implies the jump of $X(t)$
but not conversely.

Throughout this section probability $\Pv$ refers to the law of
the three processes started in distribution $\un{\hat \nu}$,
evolving according to the basic coupling as described above.

\hop

\noindent{\bf Proof of the upper bound}

\hip
 We begin with the proof of the
upper bound in Theorem \ref{tm:main}. $C$ and its variants $C_1,
C_2, \dotsc$
 denote positive constants that possibly depend on
 $\vr$ and whose values can change from line to line. We first
prove a lower bound on the density of the $\om-\eta$ second class
particles. For integers $j\in\Zb$ and $u>0$  let \be
N_j(t)=\sum_{i=j+1}^{j+2u-1}(\om_i(t)-\eta_i(t)).\label{eq:ntdef}
\ee
\begin{lm}\label{lm:2ndldp}  Let $\vr>\la>0$ and $d\geq 0$ be an
integer.  Then there are strictly positive finite
 constants $\ga=\ga(\vr)$, $C_1=C_1(\vr,d)$ and
$C_2=C_2(\vr)$  such that the following holds: if  $0<\vr-\la<\ga$,
then
 for all integers $j\in\Zb$, $u>0$  and any time $t\geq0$,
\[
\Pv\Bigl\{N_j(t)<u(\vr-\la)+d\Bigr\}\leq
C_1\exp\{-C_2u(\vr-\la)^2\}
\]
\end{lm}
\begin{proof}
For the moment, denote by $\un y(\cdot)$ a constant-rate TAZRP such
that $y_i(0)=\om_i(0)$ at all sites $i$ except for $i=0$, and $\un
z(\cdot)$ is a constant-rate TAZRP process such that
$z_i(0)=\eta_i(0)$ at all sites $i$ except for $i=0$. For $i=0$, we
pick the pair $(z_0(0),\,y_0(0))$ in distribution $\nu$,
independently of the configuration on other sites, such that $z(0)
\le y(0)\le \om^-(0)$ and $z(0) \le \eta(0)\le \om^-(0)$ holds. This
can be done since $\mu_\la\le \mu_\vr \le \hat \mu_\vr$ and
$\mu_\la\le \hat \mu_\la \le \hat\mu_\vr$. Apply the basic coupling
to ensure $\un z(t) \leq \un\eta(t) \leq\un\om(t)$ and $\un y(t)\leq
\un \om^-(t)$ for all $t\geq0$ (notice that this holds initially).
$\un y(t)$ and $\un z(t)$ are marginally time-stationary processes,
 hence we can omit the
notation for their time dependence in our arguments. However, the
pair $(\un z(t),\,\un y(t))$ is \emph{not} in product distribution
for $t>0$. Define
\[
Y=\sum_{i=j+1}^{j+2u-1}y_i\qquad\text{and}\qquad
Z=\sum_{i=j+1}^{j+2u-1}z_i,
\]
so that $N_j(t)=\sum_{i=j+1}^{j+2u-1}(\om_i(t)-y_i(t)) + Y-Z
+\sum_{i=j+1}^{j+2u-1}(z_i(t)-\eta_i(t))$. In this expression the
two sums can be estimated from above by the total number of $\om-y$
and $z-\eta$ second class particles, respectively. Using the facts
that these models only differ initially at the origin and the total
number of second class particles is preserved over time by the basic
coupling allows us to use the estimations $0\le
\sum_{i=j+1}^{j+2u-1}(\om_i(t)-y_i(t)) $ and $
\sum_{i=j+1}^{j+2u-1}(z_i(t)-\eta_i(t)) \ge z_0(0)-\eta_0(0) \ge
-\eta_0(0)$. If we plug in the lower bounds the probability of the
event increases: \[ \Pv\Bigl\{N_j(t)<u(\vr-\la)+d\Bigr\} \le
\Pv\Bigl\{Y-Z-\eta_0(0)<u(\vr-\la)+d\Bigr\} \]
Now we use exponential Markov's inequality to get for any $\al>0$
\begin{multline*}
\Pv\Bigl\{N_j(t)<u(\vr-\la)+d\Bigr\}\\
\ba
&\leq\Pv\Bigl\{\e{-\al(Y-Z-\eta_0(0))}>
\e{-\al u(\vr-\la)- d\al}\Bigr\}\\
&\leq\e{\al  u(\vr-\la)+ d\al}\cdot\Ev\bigl(\e{-\al Y}\e{\al Z}\e{\al \eta_0(0)}\bigr)\\
&\leq\e{\al  u(\vr-\la)+ d\al}\cdot\bigl[\Ev^\vr\bigl(\e{-4\al
Y}\bigr)\bigr]^{1/4}
\cdot\bigl[\Ev^\la\bigl(\e{4\al Z}\bigr)\bigr]^{1/4}\cdot\bigl[\Ev\bigl(\e{2\al \eta_0(0)}\bigr)\bigr]^{1/2}\\
&=\e{\al u(\vr-\la)+d\al}\cdot\bigl[\Ev^\vr\bigl(\e{-4\al
y_0}\bigr)\bigr]^{1/2u-1/4}
\cdot\bigl[\Ev^\la\bigl(\e{4\al z_0}\bigr)\bigr]^{1/2u-1/4}\\
&\quad \cdot\bigl[\Ev\bigl(\e{2\al \eta_0(0)}\bigr)\bigr]^{1/2}\\
&=\exp\bigl\{\al
u(\vr-\la)+d\al-\Big(\frac{1}{2}u-\frac{1}{4}\Big)\log\bigl[1+\vr(1-\e{-4\al})\bigr]
\\&\quad-\Big(\frac{1}{2}u-\frac{1}{4}\Big)\log\bigl[1+\la(1-\e{4\al})\bigr]
-\log\bigl[1+\la(1-\e{2\al})\bigr]\bigr\}\\
&\leq\exp\bigl\{-\al
u(\vr-\la)+4u(\vr+\la+\vr^2+\la^2)\al^2\\
&\quad+uC_3\al^3\! +C_4 \al + C_5
\al^2\! + \Oc(\al^3)\}. \ea
\end{multline*}
Here we used the marginal $\un\mu_\la$ and $\un\mu_\vr$ product
distributions of $\un z$ and $\un y$, and the fact that $\eta_0(0)
\sim \hat \mu_\la$ which implies $\Ev(\e{2\al \eta_0(0)} )=
(1+\la+\la\e{2\al})^{-2}$. The last inequality comes from Taylor
expansion w.r.t.\ $\al$. The $\Oc(\al^3)$ term is uniform over
$\la<\vr$ for a fixed $\vr$. In order to ensure the existence of the
exponential moment we need the condition
$\frac{\la}{1+\la}e^{2\al}<1$ to be satisfied, or equivalently, $\al
< \frac{1}{2}\log\bigl(\frac{1+\la}{\la})$. Next we pick
\[
\al=\frac{\vr-\la}{8(\vr+\la+\vr^2+\la^2)},
\]
optimizing the $\vr$ and $\la$-dependent terms in which $u$ is also
present. Comparing this to the required condition above gives a
bound for $\gamma$, i.e.\ $\gamma\le 4
\log\big(\frac{1+\vr}{\vr}\big)\cdot(\vr+\vr^2)$. By this choice of
$\al$ we get
\[\Pv\Bigl\{N_j(t)<u(\vr-\la)+d\Bigr\} \le e^{\al (d+ 2\la)+ \al^2 C_3 +\Oc(\al^3)}
\exp\big\{ -\frac{u(\vr-\la)^2}{16 (\vr +\la +\vr^2+\la^2)}\big\}\]
Now $\vr-\la \le \gamma$  gives rise to the constant $C_1$ in the
first factor and finishes the proof.
\end{proof}

Now we are ready to turn to the main estimate for the upper bound.
The idea is to bound the deviation $\Pv\{Q_a(t)\geq u+[V^\vr t]\}$
with an appropriate expression that involves the moment $\Ev\lvert
Q_a(t)-[V^\vr t]\rvert$. This is completed
 in \eqref{eq:final} below and then the upper bound comes from an
elementary integration step.

Along the way we compare currents in two processes
that we abbreviate as follows:
\[
\ba J^\vr(t)&=J^{(\frac{1}{1+\vr^2})}(t)\quad\text{for current in
the
$\un\om(.)$  process, and}\\
J^{V^\vr\!,\la}(t)&=J^{(\frac{1}{1+\vr^2})}(t)\quad\text{for current
in the $\un\eta(.)$  process.} \ea
\]
Notice that both use the same speed
$V^\vr=f'(\vr)=\frac{1}{1+\vr^2}$ for the observer.  As already
defined in the Introduction, this is the characteristic speed of the
$\un\om(.)$ process with density $\vr$.
From now on we denote by tilde the centered random variable.
\begin{lm}
 Suppose $\vr-\la<\ga$ with $\ga$
from Lemma \ref{lm:2ndldp}.  Then for positive integers $u$
and times $t\in[0,\infty)$,
\be
\ba
&\Pv\{Q_a(t)\geq2u+[V^\vr t]\}\\
&\qquad \leq \Pv\{\wt J^\vr(t)-\wt J^{V^\vr\!,\la}(t)\geq
u(\vr-\la)-t\frac{(\vr-\la)^2}{(1+\vr)^2}\}\\
&\qquad\qquad +C_1\exp\{-C_2u(\vr-\la)^2\}. \ea\label{eq:qjprob} \ee
\end{lm}
\begin{proof} Now we recall the construction of the labeled particles
and the fact that $Q_a(t)\le X(t)\quad \forall t$:
\[
\Pv\{Q_a(t)\geq2u+[V^\vr t]\}\le\Pv\{X(t)\geq2u+[V^\vr t]\}.\\
\]
$N_{[V^\vr t]}(t)$ of \eqref{eq:ntdef} counts the number of
$\om-\eta$ second class particles at time $t$ in the interval
$\{[V^\vr t]+1,\dots,[V^\vr t]+2u-1\}$. Since the second class
particles stay ordered and $X(t)$ started at the origin, the event
$\{X(t)\geq2u+[V^\vr t]\}$ implies that all these second class
particles crossed the path
 $s\mapsto [V^\vr s]+1/2$ by time $t$.
 Each such second class
particle crossing  increases $J^\vr(t)-J^{V^\vr\!,\la}(t)$ by one.
Therefore
\begin{multline*}
\Pv\{X(t)\geq2u+[V^\vr t]\}\leq\Pv\{J^\vr(t)-J^{V^\vr\!,\la}(t)
\geq N_{[V^\vr t]}(t)\}\\
\leq\Pv\Bigl\{J^\vr(t)-J^{V^\vr\!,\la}(t)\geq
u(\vr-\la)+4\vr+1\Bigr\}+\Pv\Bigl \{N_{[V^\vr
t]}(t)<u(\vr-\la)+4\vr+1\Bigr\}.
\end{multline*}
Combine the previous displays to get
\begin{align}
&\Pv\{Q_a(t)\geq2u+[V^\vr t]\}\leq
\Pv\Bigl\{J^\vr(t)-J^{V^\vr\!,\la}(t)
\geq u(\vr-\la)+4\vr+1\Bigr\}\label{eq:temp7}\\
&\qquad+\Pv\Bigl\{N_{[V^\vr
t]}(t)<u(\vr-\la)+4\vr+1\Bigr\}\label{eq:nupr}
\end{align}
To line \eqref{eq:nupr} apply Lemma \ref{lm:2ndldp} with
$d=\lceil4\vr+1\rceil$. We see that line  \eqref{eq:nupr} is bounded
by the last exponential term in \eqref{eq:qjprob}. Recall the
definition of the flux \eqref{eq:flux}.
 If $\un\om$ and $\un\eta$ start from
their respective $\un\mu_\vr$ and $\un\mu_\la$ equilibria then we
would have, due to \eqref{eq:jvt} and \eqref{eq:Qmean},
\[
\begin{split}
\Ev^\vr(J^\vr(t))-\Ev^\la(J^{V^\vr\!,\la}(t))
&=t\bigl( f(\vr)-f(\la)-f'(\vr)(\vr-\la)\bigr)\\
&=t\frac{(\vr-\la)^2}{(1+\vr)^2(1+\la)}
\end{split}
\]
where we ignored the error coming from the integer part of $V^\vr
t$. This error can be estimated from above by the number of
particles sitting at one site in each of the processes, i.e.\
$\Ev^\vr(J^\vr(t))-\Ev^\la(J^{V^\vr\!,\la}(t))-t\frac{(\vr-\la)^2}{(1+\vr)^2(1+\la)}\le
\Ev^\vr(\om_i)+\Ev^\la(\eta_i)=\vr+\la$. Our processes are also
perturbed initially at the origin by $\hat \nu$ being different of
$\nu$, which gives an error $\Ev (J^\vr(t)) - \Ev^\vr(J^\vr(t)) \le
\Ev(\om_0(0))= 2\vr+1$. In the $\un\eta$ process the term $-
\Ev(J^{V_\vr,\la}(t))+\Ev^\la(J^{V_\vr,\la}(t))$ is negative so it
can be estimated from above by $0$. The term $4\vr+1$ inside the
probability on line \eqref{eq:temp7}
 makes up for these errors. So, we get $ \Pv\{\wt J^\vr(t)-\wt J^{V^\vr\!,\la}(t)\geq
u(\vr-\la)-t\frac{(\vr-\la)^2}{(1+\vr)^2(1+\la)}\}$
 Finally, we omit the factor $1+\la$ in
 the denominator, implying the decrease of the right hand side of
 the inequality, and increasing the probability.
\end{proof}

\begin{lm}\label{lm:vardif}
\[
\ba \Vv(J^{\vr}(t))&\leq 12\vr(\vr+1) +2+
2\vr(1+\vr)\Ev\bigl(|[V^\vr
t]-Q_a(t)|\bigr),\\
 \Vv(J^{V^\vr\!,\la}(t))&\leq 12\vr(\vr+1) + 4t(\vr-\la)+2\vr(1+\vr)
\Ev\bigl(|[V^\vr t]-Q_a(t)|\bigr).\ea
\]
\end{lm}
\begin{proof}
The variance $\Vv$
 in the statement is taken in the three-process coupling where the geometric distribution  $\un\mu_\vr$
 is initially perturbed at the origin.
 Recall that $\Vv^\vr$ denotes
variance in the stationary process $\un y(t)$ coupled to $\un
\om(t)$, and denote by $J_y^\vr(t)$ the height function in this
process at site $[V^\vr t]$ at time $t$. The following estimation
holds: \be
\begin{split}
\Vv(J^\vr(t))&\le 2 \Vv(J^\vr(t)-J_y^\vr(t))+ 2 \Vv^\vr J^\vr(t) \\
&\le 2\Ev(\om_0(0)^2) + 2 \Vv^\vr(\om_0(0))\cdot \Ev(|Q_a(t)-[V^\vr t]|)\\
&\le 12\vr(\vr+1)+ 2  + 2\vr(\vr+1)\Ev(|Q_a(t)-[V^\vr t]|).
\end{split}\label{eq:temp12}\ee
Here we applied Proposition \ref{pr:se} to the second term
$\Vv^\vr(J^\vr(t))$ in the first line and used the fact that
$J^\vr(t)$ and $J_y^\vr(t)$ only differ at the origin by
$\om_0(0)-y_0(0)$, so its variance can be estimated from above by
the second moment of $\om_0(0)$, equal to the second moment of $1$
plus the sum of two independent geometrically distributed random
variables with mean $\vr$. Similar computation for the $\la$-density
process is as follows (with stationary $\la$-density process $\un
z(t)$ and its height function $J_z^{V^\vr,\la}(t)$):
\[
\begin{split}
\Vv(J^{V^\vr\!,\la}(t)) &\le 2 \Vv(J^{V^\vr\!,\la}(t)-J_z^{V^\vr,\la}(t))+ 2 \Vv^\la J^{V^\vr\!,\la}(t)\\
&\le2\Ev(\eta_0(0)^2) + 2 \Vv^\la(\eta_0(0))\cdot \Ev(|Q^\la(t)-[V^\vr t]|)\\
&\le 4\la(\la+1) + 8 \la^2+ 2\la(\la+1)\Ev(|Q^\la(t)-[V^\vr t]|)\\
&\le 12\vr(\vr+1)+
2\la(\la+1)\Ev(|Q^\la(t)-Q_a(t)|)\\
&\quad \quad +2\vr(\vr+1)\Ev(|Q_a(t)-[V^\vr
t]|)\\
&\le12\vr(\vr+1) +4(\vr-\la)t+2\vr(\vr+1)\Ev(|Q_a(t)-[V^\vr t]|).
\end{split}
\]
In the first line we cut the variance into two terms, where the
first term can be estimated from above by $ 2\Ev(\eta_0(0)^2)$, and
for the second term we applied Proposition \ref{pr:se} with a second
class particle $Q^\la(t)$ added to the $\un \eta$ process. Then we
used triangle inequality and Proposition \ref{pr:se} again to get
$\Ev(|Q^\la(t)-Q_a(t)|)=\Ev(Q^\la(t)-Q_a(t))=t\frac{1}{(1+\la)^2}-t\frac{1}{(1+\vr)^2}$.
 The absolute value can be omitted by $Q^\la(t)\ge Q_a(t)$, which we show below.
 Initially $Q^\la(0)=Q_a(0)=0$.
 The jump rates of $Q^\la(t)$ and
$Q_a(t)$ are $\mathbb I_{\eta_{Q^\la(t)}(t)=0}$ and $\mathbb
I_{\om^-_{Q_a(t)}(t)=0}$, respectively. Second class particles are
only allowed to jump to the right, and a right jump of $Q_a$ without
$Q^\la$ is impossible by the basic coupling and the above rates
since $\eta_{i}(t)\leq\om^-_{i}(t)$.

Some algebra then leads to the estimation $2\la (\la+1)
\Ev(|Q^\la(t)-Q_a(t)|) \le 4(\vr-\la)t$. Collecting terms completes
the proof of the lemma.
\end{proof}

We come to the lemma that summarizes all the previous estimations.

\begin{lm}
For any real $u\geq 1$ and time $t> 0$,
\begin{multline}
\Pv\{Q_a(t)\geq2u+[V^\vr t]\}\leq C_3\frac{t^2}{u^4}\cdot\Ev\bigl(|Q_a(t)-[V^\vr t]|\bigr)\\
+C_4\frac{t^2}{u^3}+C_5\frac{t^2}{u^4}
+C_1\exp\Bigl\{-C_2\frac{u^3}{t^2}\Bigr\}+\e{-u}.\label{eq:final}
\end{multline}
\end{lm}
\begin{proof}
Set $b=\frac{2}{(1+\vr)^2}(\gamma\wedge\vr)$ where $\gamma$ is the
constant from Lemma \ref{lm:2ndldp}.  We proceed by considering
three cases.

{\sl Case 1:}  $1\le u<bt$.  Suppose first $u$ is an integer as it
was in the proof of \eqref{eq:qjprob}.  Throughout density $\vr$ has
been fixed, and now we also fix
\[
\la=\vr-\frac{u}{2t}(1+\vr)^2.
\]
The constraint on $u$ guarantees that $\la>0$ and
$\vr-\la<\ga$ which was required for \eqref{eq:qjprob}. The point
of this choice of $\la$ is to maximize
 the lower bound inside the
probability on the right-hand side of \eqref{eq:qjprob}.
So, continuing from \eqref{eq:qjprob} with
 Chebyshev's inequality and the previous lemma,
\begin{multline*}
\Pv\{Q_a(t)\geq2u+[V^\vr t]\}\\
\ba &\leq \Pv\Bigl\{\wt J^\vr(t)-\wt J^{V^\vr\!,\la}(t)\geq
\frac{u^2}{4t}(1+\vr)^2\Bigr\}+
C_1\exp\Bigl\{-C_2\frac{u^3}{t^2}\Bigr\}= \\
&\leq C_3\frac{t^2}{u^4} \Vv(J^\vr(t)-J^{V^\vr\!,\la}(t))+
C_1\exp\Bigl\{-C_2\frac{u^3}{t^2}\Bigr\}\\
&\leq C_3\frac{t^2}{u^4}(2\Vv(J^\vr(t)+2 \Vv J^{V^\vr\!,\la}(t))+
C_1\exp\Bigl\{-C_2\frac{u^3}{t^2}\Bigr\}\\
&\leq C_3\frac{t^2}{u^4}\cdot\Ev\bigl(|Q_a(t)-[V^\vr t]|\bigr)+
C_4\frac{t^2}{u^3}+C_5\frac{t^2}{u^4}+C_1\exp\Bigl\{-C_2\frac{u^3}{t^2}\Bigr\}.
\ea
\end{multline*}
Extension  from integral $u$ to real $u$ is achieved by adjusting
constants on the last line above.

{\sl Case 2:}  $bt\le u< 2t$. Then $bu/2< bt$ and
\[
\Pv\{Q_a(t)\geq2u+[V^\vr t]\}\leq\Pv\{Q_a(t)\geq2\cdot bu/2+[V^\vr
t]\}.
\]
{\sl Case 1}  can be applied
with $u$ replaced by $bu/2$, at the price of adjusting
some constants with powers of $b$.

{\sl Case 3:}  $u\ge 2t$.  Since  $Q_a(t)$ is bounded above by a
rate one Poisson process $N(t)$,
\[
\begin{split}
 \Pv\{Q_a(t)\geq2u+[V^\vr
t]\}&\leq\Pv\{Q_a(t)\geq2u\}\leq \Pv\{N(t)\geq 2u\}\\
&\le e^{-2u}\cdot \Ev(e^{N(t)})=e^{t(e-1)-2u} \le e^{-u}
\end{split}
\]
for all $t$.  Combining the bounds from the cases proves the
lemma.
\end{proof}

We are ready to complete the proof of the upper bound of Theorem
\ref{tm:main}. Introduce the notation \[
\Psi(t)=\Ev\bigl(\,|Q_a(t)-[V^\vr t]|\,\bigr).\]

In order to get a bound on $\Pv\{|Q_a(t)-[V^\vr t]|\ge 2u\}$ we need
a lower tail bound for $Q_a(t)$, similar to \eqref{eq:final}. This
can be obtained by methods analogous to the ones we have applied
throughout Section \ref{sc:ub}. In this case we have to couple a
process with initial distribution $\un{\hat \mu}_\sigma$ for some
$\sigma>\vr$ close enough to $\vr$ and derive similar lemmas to get
the same conclusion for $\Pv\big\{Q_a(t) \le -2u +[V^\vr t]\big\}$
as in (\ref{eq:final}) with readjusted constants.

Introduce  a large constant  $2<\kon<\infty$. For $t\ge 1$ and $u\ge
\kon t^{2/3}$ we can combine \eqref{eq:final} and the matching lower
tail bound. Replace $2u$ by $u$ (we made sure $u/2\ge 1$), then some
algebra leads to \be
\begin{split}
&\Pv\{\,\lvert Q_a(t)-[V^\vr t] \rvert \ge u\} \\
&\qquad \leq C_1\frac{t^2}{u^4} \Psi(t)
+C_2(\kon)\Bigl(\,\frac{t^2}{u^3} +
\exp\Bigl\{-C_3(\kon)\frac{u}{t^{2/3}}\Bigr\} \,\Bigr).
\end{split}\label{eq:final2}\ee
The two exponential terms in \eqref{eq:final} were combined via
$e^{-u}$ $\le$ $\exp(-ut^{-2/3})$ and $\exp(-C_2u^3t^{-2})\le
\exp(-C_2\kon^{2}ut^{-2/3})$.

Let $1\le \mom<3$.  Integrate bound \eqref{eq:final2} over
$u\in [\kon t^{2/3},\infty)$:
\be\begin{split}
&\Ev\bigl(\,|Q_a(t)-[V^\vr t]|^\mom\,\bigr)\\
&\qquad \le \kon^\mom t^{2\mom/3} +\mom\int_{\kon t^{2/3}}^\infty
\Pv\{\,\lvert Q_a(t)-[V^\vr t] \rvert \ge u\} u^{\mom-1} \,\di u\\
&\qquad \le C_1 \kon^{\mom-4} \Psi(t) t^{2+(2/3)(\mom-4)} \;+\;
C_4(\kon)t^{2\mom/3}.
\end{split}\label{eq:final3}\ee
The constant $C_1$ depends on $\mom$ but not on $\kon$, while $C_4$
depends on both $\mom$ and $\kon$. To get the final bounds, take
first $\mom=1$ in \eqref{eq:final3} to get
\[
\Psi(t) \le C_1 \kon^{-3} \Psi(t)  \;+\; C_4(\kon)t^{2/3}.
\]
Since $C_1$ is independent of $\kon$, fixing  $\kon$ large enough
gives $\Psi(t)\le C_5(\kon)t^{2/3}$. Using this bound  in the last
line of \eqref{eq:final3} implies
\[
\Ev\bigl(\,|Q_a(t)-[V^\vr t]|^\mom\,\bigr) \le C_6(\kon)t^{2\mom/3}
\]
for $1<\mom<3$, which completes the proof of the upper bound of
Theorem \ref{tm:main}.

\section{Lower bound}

The lower bound is proved by perturbing the stationary distribution
of a process on a segment of the lattice. We start with discussing
the initial distribution of some coupled processes.


\hop
{\bf Perturbing a segment initially}

\hip Recall again the characteristic speeds
$V^\vr=\frac{1}{(1+\vr)^2}$ and $V^\la=\frac{1}{(1+\la)^2}$. We set
$\vr>\la$, hence $V^\vr<V^\la$. Throughout this section $u>0$
denotes a fixed positive integer, and
\[
n=[V^\la t]-[V^\vr t]+u.
\]
Recall the Definition \ref{def:measures} of the measures $\nu$ and
$\hat\nu$. Define an initial product distribution of two
configurations $(\un\eta(0),\un\ze(0))$
 by setting the marginals
on each lattice site $i$:
\[\begin{cases}
(\eta_i(0),\,\ze_i(0))\sim\nu &\text{if }i<-n,\\
(\eta_i(0),\ \ze_i(0))\sim\hat \nu &\text{if }i=-n,\\
\eta_i(0)=\ze_i(0)\sim\mu_\la&\text{if }-n<i\leq0,\\
(\eta_i(0),\,\ze_i(0))\sim\nu &\text{if }i>0.
\end{cases}
\]
Note that the number of particles at site $-n$ is $\hat\mu_\la$-
distributed for $\un\eta(0)$. Except for this perturbation
$\un\eta(0)$ starts in the stationary Geometric product distribution
$\un \mu_\la$. The process $\un\ze(\cdot)$ initially has
distribution $\un\mu_\vr$, except at sites $\{-n+1,\,\dots,\,0\}$
where the parameter $\vr$ has been replaced by $\la$, and at site
$-n$ where it has measure $\hat\mu_\vr$.

We add a second class particle to the process $\un \eta(.)$
initially started at site $-n$ and denote its position at time $t$
by $Q^{(-n)}(t)$. Introduce, as before, $\un \eta^+(t):= \un\eta(t)+
\un \delta_{Q^{-n}(t)}$.

Define a third initial configuration by
\[
\xi_i(0)=\left\{
\ba
&\ze_i(0)&&\text{if }i\leq-n,\\
&\eta_i(0)&&\text{if }i>-n.
\ea
\right.
\]
Apply the basic coupling to obtain the joint evolution of all these
processes. This guarantees the majorizations
\[
\un\eta(t)\leq\un\xi(t)\leq\un\ze(t)\qquad\text{and}\qquad \un h^\ze(t)\leq\un h^\xi(t)
\]
where the last inequality is for column heights.

As before in \eqref{eq:jvdef}, denote the net particle currents by
$J^{V,\eta}$ and $J^{V,\ze}$ in the respective processes
$\un\eta(\cdot)$ and $\un\ze(\cdot)$. The first observation is that
 $Q^{(-n)}$ gives one-sided control over the difference
of  these currents.
\begin{lm}\label{lm:geom}
For any $V\in\Rb$
\[
Q^{(-n)}(t)\leq[Vt] \quad\text{implies}\quad
J^{V,\ze}(t)-J^{V,\eta}(t)\leq 0.
\]
\end{lm}
\begin{proof}
Denote the positions of the  $\xi-\eta$ second class particles at time
$t$  by
\[ \dotsm \le Y_k(t)\le\dotsm \le Y_{-2}(t)\le Y_{-1}(t)\le Y_0(t).\]
This order of the labels is again preserved over time. Initially
$Y_0(0)\le -n=Q^{(-n)}(0)$, and the respective jump  rates $\mathbb
I_{\eta_{Y_0(t)}=0}$ and $\mathbb I_{\eta_{Q^{(-n)}(t)}=0}$ of the
second class particles ensure $Y_0(t)\le Q^{(-n)}(t)$ for all later
times. Indeed, from the time $Y_0$ has jumped on $Q^{(-n)}$ they
always jump together. So, $Y_0$ cannot leave $Q^{(-n)}$ behind
itself. Saying it another way, once $Q^{(-n)}(t)=Y_0(t)$, this
property will hold forever.

Now note that there are no $\xi-\eta$ second class particles
strictly to the right of $Q^{(-n)}(t)$ at time $t$. So, the height
difference is zero to the right of $Q^{(-n)}(t)$. Also recall the
majorization $\un h^\ze(t)\leq\un h^\xi(t)$. Thus
 under $\{Q^{(-n)}(t)\leq[Vt]\}$ we have
\[
0=h^\xi_{[Vt]}(t)-h^\eta_{[Vt]}(t) \geq
h^\ze_{[Vt]}(t)-h^\eta_{[Vt]}(t)=J^{V,\ze}(t)-J^{V,\eta}(t).
\]
\end{proof}

Let $\un{\hat\om}(\cdot)$ be a TAZRP started from the distribution
$\un{\hat\mu}_\vr$ shifted $n$ sites to the left, i.e.\  at site
$-n$ the number of particles is $\hat \mu_\vr$-distributed. The next
lemma compares the distributions of $\un\ze$ and $\un{\hat\om}$.
\begin{lm}\label{lm:rn}
Denote by $\Pv^{\hat\om}$ and $\Pv^\ze$ the probability of events
that depend only on the respective processes $\un{\hat\om}(\cdot)$
and $\un\ze(\cdot)$. Then there exist $\gamma>0$ such that for all
$\la>\vr-\gamma$ the following inequality holds:
\[
\Pv^\ze(\cdot)\leq\Pv^{\hat\om}(\cdot)^\frac12\cdot\exp\Bigl[\frac{n(\vr-\la)^2}
{\vr(1+\vr)}\Bigr].
\]
\end{lm}
\begin{proof}
Let
\[
Z=\sum_{i=-n+1}^0\ze_i(0).
\]
 $Z$ has a Negative Binomial($n, \frac{1}{1+\la}$) distribution,
 namely
 \[m^\la(z):=\Pv(Z=z)=\binom{n+z-1}{n-1}\Big(\frac{\la}{1+\la}\Big)^z
 \Big(\frac{1}{1-\la}\Big)^{n}\]
for $z\ge 0$. We use the Cauchy-Schwarz inequality below to perform
a change of measure on this negative binomial distribution.
\[
\ba \Pv^\ze(\cdot)&=\sum_{z=0}^\infty\Pv^\ze(\,\cdot\,|\,Z=z)\,
[m^\vr(z)]^\frac12\cdot\frac{m^\la(z)}{[m^\vr(z)]^\frac12}\\
&\leq\Bigl[\sum_{z=0}^\infty[\Pv^\ze(\,\cdot\,|\,Z=z)]^2\,
m^\vr(z)\Bigr]^\frac12\cdot\Bigl[\sum_{z=0}^\infty\frac{[m^\la(z)]^2}
{m^\vr(z)}\Bigr]^\frac12. \ea
\]
Now the last part equals $\sum_{z=0}^\infty \binom{n+z-1}{n-1}\bigl(
\frac{\la^2(1+\vr)}{(1+\la)^2\vr}\bigr)^z\cdot\bigl(
\frac{(1+\vr)}{(1+\la)^2}\bigr)^n$. In this formula we can recognize
the negative binomial mass function with parameter
\\$0<1-\frac{\la^2(\vr+1)}{(1+\la)^2 \vr}< 1$, so we use the identity \[\sum_{z=0}^\infty \binom{n+z-1}{n-1}\Bigl(
\frac{\la^2(1+\vr)}{(1+\la)^2\vr}\Bigr)^z\cdot \Bigl(1-
\frac{\la^2(1+\vr)}{(1+\la)^2\vr}\Bigr)^n=1\] and the remaining
factors equal $\Big(1-\frac{\la^2(\vr+1)}{(1+\la)^2
\vr}\Big)^{-n}\cdot
\Big(\frac{\vr+1}{(\la+1)^2}\Big)^n=\Bigl[1-\frac{(\vr-\la)^2}{\vr(1+\vr)}\Bigr]^{-n}$.
So we get
\[
\ba
\Pv^\ze(\cdot)&\le\Bigl[\sum_{z=0}^\infty[\Pv^\ze(\,\cdot\,|\,Z=z)]^2\,m^\vr(z)\Bigr]^\frac12
\cdot\Bigl[1-\frac{(\vr-\la)^2}{\vr(1+\vr)}\Bigr]^\frac{-n}{2}\\
&\leq\Bigl[\sum_{z=0}^n\Pv^\ze(\,\cdot\,|\,Z=z)\,m^\vr(z)\Bigr]^\frac12
\cdot\exp\Bigl[\frac{n(\vr-\la)^2}{\vr(1+\vr)}\Bigr]. \ea
\]
In the last line we used the fact that $1-x \ge e^{-2x}$ for all
$0<x\le0.5$. This will be satisfied for
$x=\frac{(\vr-\la)^2}{\vr(\vr+1)}$, for $\la$ close enough to $\vr$.
All that is left is to recognize that $\Pv^\ze(\,\cdot\,|\,Z=z)$ is
the probability depending on a process $\un\ze(\cdot)$ whose initial
distribution coincides with the distribution of $\hat \om$ outside
$\{-n+1\dots0\}$, with $z$ particles distributed in that interval
with each configuration equally likely. Namely, each configuration
with $z$ particles in this interval has probability
$\Big(\frac{\la}{1+\la}\Big)^z \Big(\frac{1}{1+\la}\Big)^n$, so
conditioned on the event that $z$ particles are here initially, each
configuration has probability $\frac{1}{\binom{n+z-1}{n-1}}$.
Summing these conditionals, weighted with the Negative
Binomial$(n,\,\frac{1}{1+\vr})$ mass function $m^\vr(z)$, gives the
product Geometric initial distribution $\un\mu_\vr$ of the
 process $\un{\hat\om}(\cdot)$.
\end{proof}


\hop
{\bf Proof of the lower bound}

\hip Now we only need two steps to obtain the lower bound. In this
part we need again the concept of a second class antiparticle
started from the origin on a $\un\mu_\vr$-stationary process
$\un\om(\cdot)$
 initially perturbed by setting $\om_0(0) \sim
\hat\mu_\vr+\delta_0$. We denote its position at time $t$ by
$Q_a(t)$. The quantity of primary interest is abbreviated, as
before, by $\Psi(t)=\Ev(|Q_a(t)-[V^\vr t]|)$.

The first step is proving an upper bound on $Q^{(-n)}(t)$. As before
$u$ is an arbitrary but fixed positive integer
 and $n=[V^\la t]-[V^\vr t]+u$.
\begin{lm}
\be\label{eq:qgevr} \Pv\{Q^{(-n)}(t)>[V^\vr
t]\}\leq\frac{\Psi(t)}{u}+\frac{4t(\vr-\la)}{u}+\frac{2}{u}. \ee
\label{lm:LBlm3} \end{lm}
\begin{proof}
Below $Q^\la(t)$ stands for the position of a second class particle
started from the origin on a process $\un{\hat\eta}(\cdot)$ in
$\un{\hat\mu}_\la$ distribution. Translation invariance implies
\begin{multline*}
\Pv\{Q^{(-n)}(t)>[V^\vr t]\}\\
\ba &=\Pv\{Q^{(-n)}(t)+n-[V^\la t]>u\}\\&=\Pv\{Q^\la(t)-[V^\la
t]>u\}\\&\leq\frac{\Ev(|Q^\la(t)-[V^\la t]|)}{u}\\
&\leq\frac{\Ev(|Q^\la(t)-Q_a(t)|)}{u}+\frac{\Ev(|Q_a(t)-[V^\vr
t]|)}{u}+\frac{[V^\la t]-[V^\vr t]}{u}. \ea \end{multline*}  As in
Lemma \ref{lm:vardif}, the first term equals
$\frac{t}{u}\big(\frac{1}{(1+\la)^2}-\frac{1}{(1+\vr)^2}\big)$,
bounded from above by $\frac{2}{u}t(\vr-\la)$ after some
calculations. The second term is $\Psi(t)/u$, and the third term is
similarly estimated by $\frac{2}{u}t(\vr-\la)+\frac{2}{u}$, the last
part coming from possible integer part errors.
\end{proof}

The second step is an estimate of the probability of the complement
of the event in \eqref{eq:qgevr}.
\begin{lm}
For any $0<K<t\frac{(\vr-\la)^2}{(1+\vr)^3}-2\vr$,
\[
\ba \Pv\{Q^{(-n)}(t)\leq[V^\vr
t]\}&\leq\frac{\big(2\vr(1+\vr)\big)^{1/2}\big(6+\Psi(t)\big)^{1/2}}
{t\frac{(\vr-\la)^2}{(1+\vr)^3}-2\vr-K}\cdot\exp\Bigl[\frac{n(\vr-\la)^2}{\vr(1+\vr)}\Bigr]\\
&\quad+\frac{2\vr(1+\vr) \big(6+\Psi(t)\big)}{K^2}+
\frac{4t(\vr-\la)}{K^2}. \ea
\]
\label{lm:LBlm4} \end{lm}
\begin{proof}
Lemma \ref{lm:geom} leads to
\begin{align}
\Pv\{Q^{(-n)}(t)\leq[V^\vr
t]\}&\leq\Pv\{J^{V^\vr\!,\ze}(t)-J^{V^\vr\!,\eta}(t)
\leq 0\}\notag\\
&\leq\Pv\{J^{V^\vr\!,\ze}(t)\leq K+t\big(\frac{\la}{1+\la}-\frac{\la}{(1+\vr)^2}\big)+2\la\}\label{eq:torn}\\
&\quad
+\Pv\Bigl\{J^{V^\vr\!,\eta}(t)>K+t\big(\frac{\la}{1+\la}-\frac{\la}{(1+\vr)^2}\big)+2\la
\Bigr\} \label{eq:tovarj}
\end{align}
We apply Lemma \ref{lm:rn} to line \eqref{eq:torn} to bound it by
the pro\-ba\-bi\-li\-ty of the process $\un{\hat\om}$:
\[
\ba (\ref{eq:torn})&\le \bigl[\Pv^{\hat\om}\{J^\vr(t)\leq
K+t\big(\frac{\la}{1+\la}-\frac{\la}{(1+\vr)^2}\big)+2\la\}
\bigr]^\frac12
\cdot\exp\Bigl[\frac{n(\vr-\la)^2}{\vr(1+\vr)}\Bigr]\\
&\leq\bigl[\Pv^{\hat\om}\{\wt J^\vr(t)\leq
K-t\frac{(\vr-\la)^2}{(1+\la)(1+\vr)^2}+2\la\}\bigr]^\frac12
\cdot\exp\Bigl[\frac{n(\vr-\la)^2}{\vr(1+\vr)}\Bigr]\\
&\leq\bigl[\Pv^{\hat\om}\{\wt J^\vr(t)\leq
K-t\frac{(\vr-\la)^2}{(1+\vr)^3}+2\vr\}\bigr]^\frac12
\cdot\exp\Bigl[\frac{n(\vr-\la)^2}{\vr(1+\vr)}\Bigr]\\
&\leq\frac{[\Vv(J^\vr(t))]^{1/2}}{t\frac{(\vr-\la)^2}{(1+\vr)^3}-2\vr-K}
\cdot\exp\Bigl[\frac{n(\vr-\la)^2}{\vr(1+\vr)}\Bigr]\\
&\leq\Big\{\frac{\big(2\vr(1+\vr)\big)^{1/2}\big(6+\Psi(t)\big)^{1/2}}{t\frac{(\vr-\la)^2}{(1+\vr)^3}-2\vr-K}
+\frac{\sqrt{2}}{t\frac{(\vr-\la)^2}{(1+\vr)^3}-2\vr-K}\Big\}
\cdot\exp\Bigl[\frac{n(\vr-\la)^2}{\vr(1+\vr)}\Bigr]. \ea
\]
When centering, we used $J^\vr(t)=\wt J^\vr(t)+[ \Ev J^\vr(t) -
\Ev^\vr J^\vr(t)]+ [\Ev^\vr J^\vr(t) -
t\big(\frac{\vr}{1+\vr}-\frac{\vr}{(1+\vr)^2}\big)]+t\big(\frac{\vr}{1+\vr}-\frac{\vr}{(1+\vr)^2}\big)$,
where both error terms in the brackets $[\cdot]$ are nonnegative.
So, $J^\vr(t)\leq
K+t\big(\frac{\la}{1+\la}-\frac{\la}{(1+\vr)^2}\big)+2\la$ implies
$\wt J^\vr(t)\leq K-t\frac{(\vr-\la)^2}{(1+\la)(1+\vr)^2}+2\la$.
Then, in the last line we used Lemma \ref{lm:vardif} to bound the
variance by the function $\Psi(t)$ even though the second class
particle in $\un{\hat\om}$ starts
 at $-n$ rather than at the origin, namely $\Vv(J^\vr(t))\le 2+12\vr(1+\vr)+2\vr(1+\vr)\Psi(t)$.
 We use the inequality
$\sqrt{a+b}\le \sqrt{a}+\sqrt{b}$ in order to get rid of the term
$+2$. The only change needed in the proof of Lemma \ref{lm:vardif}
is in the calculation \eqref{eq:temp12} where we must add
$\Ev(\om_{-n}(0)^2)$ instead of $\Ev(\om_{0}(0)^2)$, but these are
equal.

 A similar coupling consideration
shows that $ \Ev(J^{V^\vr\!,\eta}(t))$ differs by at most $\la$ from
the same expectation taken under a stationary $\un\mu_\la$ initial
condition. Thus taking integer parts again into account, giving
another error term $\la$, line \eqref{eq:tovarj} is bounded from
above by
\[
\Pv\Bigl\{\wt J^{V^\vr\!,\eta}(t)>K\Bigr\}\leq
\frac{\Vv(J^{V^\vr\!,\eta})}{K^2}.
\]
 Lemma \ref{lm:vardif} can be applied again to bound this variance
 with the similar change of shifting all processes by $n$ sites
 to the left.
Hence we can continue from above to bound
 line \eqref{eq:tovarj} with
\[
\frac{2\vr(1+\vr) \big(6+\Psi(t)\big)}{K^2}+\frac{4t(\vr-\la)}{K^2}.
\]

\end{proof}

Now we are at the last step of proving  the
 lower bound of Theorem \ref{tm:main}.  By Jensen's inequality it suffices to prove for the case $\mom=1$,
in other words that
\[
\liminf_{t\to\infty} t^{-2/3}\Psi(t)>0.
\]
In the last two lemmas take
\[
u=\lc ht^{2/3}\rc,\quad\vr-\la=bt^{-1/3}, \quad \text{ and }\quad K=bt^{1/3},
\]
where $h$ and  $b$ are large, in particular $b$ large enough to have
$bt^{\frac{1}{3}}<\frac{b^2}{(1+\vr)^3} t^{\frac{1}{3}}-2\vr$ so
that $K$ satisfies the assumption of Lemma \ref{lm:LBlm4}. Then
\[n=[V^\la t]-[V^\vr t]+u\le
\frac{2+\vr+\la}{(1+\la)^2(1+\vr)^2}(\vr-\la)t+2+u\le (2b+h) t^{2/3}
\le C t^{2/3} \] for large enough $t$. We can simplify the outcomes
of Lemma \ref{lm:LBlm3} and Lemma \ref{lm:LBlm4}  to the
inequalities \be \Pv\{Q^{(-n)}(t)>[V^\vr t]\}\leq
C\frac{\Psi(t)}{t^{2/3}} +\frac{4b}{h} +\frac{2}{h t^{2/3}}
\label{eq:LBaux3}\ee and \be\begin{split} \Pv\{Q^{(-n)}(t)\leq[V^\vr
t]\}&\leq C\biggl(\frac{6+\Psi(t)}{t^{2/3}}\biggr)^{1/2}
+C\frac{6+\Psi(t)}{t^{2/3}}+ \frac{4}{b}+ \frac{C}{t^{1/3}}.
\end{split}
\label{eq:LBaux4}\ee
The new constant $C$ depends on $b$ and $h$.

The lower bound now follows because the left-hand sides of
\eqref{eq:LBaux3}--\eqref{eq:LBaux4} add up to one for each
fixed $t$, while we can fix
 $b$ large enough and then $h$ large enough so
that $4b/h + 4/b < 1$.  Then $t^{-2/3}\Psi(t)$ must have a positive
lower bound for all large enough $t$. This completes the proof of
Theorem \ref{tm:main}.

\section*{Acknowledgments}

The authors wish to thank Timo Sepp\"al\"ainen for fruitful discussions and advice on the subject.

\hop
\begin{center}
\footnotesize{
\begin{tabular}{l}
{\sc M\'arton Bal\'azs}\\
{\sc MTA-BME Stochastics Research Group, Institute of Mathematics}\\
{\sc Budapest University of Technology and Economics}\\
{\sc 1 Egry J\'ozsef u., $5^{\text{th}}$ floor 7, Bld.\ H}\\
{\sc 1111 Budapest, Hungary}\\
{\tt balazs@math.bme.hu}\\
\\
{\sc J\'ulia Komj\'athy}\\
{\sc BME, Institute of Mathematics}\\
{\sc Budapest University of Technology and Economics}\\
{\sc 1 Egry J\'ozsef u., $5^{\text{th}}$ floor 7, Bld.\ H}\\
{\tt komyju@math.bme.hu}
\end{tabular}
}
\end{center}

\bibliography{refsmarton}
\bibliographystyle{plain}

\end{document}